\documentclass{amsart}
\pdfoutput=1
\usepackage{amsthm, amssymb, array, enumerate, times, txfonts, mathrsfs, bbm, verbatim}
\usepackage[centredisplay, PostScript=dvips]{diagrams}
\usepackage[initials]{amsrefs}
\title{Enumerating the Derangements of an $n$-Cube via M\"obius 
Inversion}
\author{Colin G.Bailey}
\address{School of Mathematics,  Statistics \& Operations Research\\
Victoria University of Wellington\\
PO Box 600\\
Wellington\\
NEW ZEALAND
}
\email{Colin.Bailey@vuw.ac.nz}
\author{Joseph S.Oliveira}
\address{
Pacific Northwest National Laboratories\\
Richland\\
U.S.A.}
\email{Joseph.Oliveira@pnl.gov}
\date{2009, February 3}

\let\rsf\mathscr

\def\one{\mathbf1}
\def\zero{\mathbf0}
\def\Fr{\operatorname{Fr}}
\def\Fix{\operatorname{Stab}}
\def\Aut{\operatorname{Aut}}
\def\cork{\operatornamewithlimits{co-rk}}
\def\rk{\operatornamewithlimits{rk}}
\def\dim{\operatornamewithlimits{dim}}
\def\stab{\operatornamewithlimits{stab}}
\def\id{\operatorname{id}}
\def\At{\operatorname{At}}
\def\caret{\mathbin{\hat{\hphantom{m}}}}
\def\uparr{\!\uparrow}
\providecommand{\meet}{\mathbin{\wedge}}
\providecommand{\join}{\mathbin{\vee}}
\newcommand{\card}[1]{\left| #1\right|}
\newcommand{\comp}[1]{\overline{#1}}
\ifx\upharpoonright\undefined
     \def\restrict{\hbox{\rm\kern0.166em\accent"12\kern-0.536em$\vert$\kern0.3em}}%
  \else
     \def\restrict{\upharpoonright}%
  \fi
\makeatletter
\def\rightharpoonupfill@#1{\setboxz@h{$#1-\m@th$}\ht\z@\z@
  $#1\m@th\copy\z@\mkern-6mu\cleaders
  \hbox{$#1\mkern-2mu\box\z@\mkern-2mu$}\hfill
  \mkern-6mu\mathord\rightharpoonup$}
\def\overrightharpoonup{\mathpalette\overrightharpoonup@}
\def\overrightharpoonup@#1#2{\vbox{\ialign{##\crcr\rightharpoonupfill@#1\crcr
 \noalign{\kern-\ex@\nointerlineskip}$\m@th\hfil#1#2\hfil$\crcr}}}

\def\twoSet#1#2{\left\{%
\vphantom{#2}#1\thinspace\right|\nolinebreak[3]\left.%
  #2%
  \vphantom{#1}%
  \right\}%
}
\def\oneSet#1{\left\lbrace#1\right\rbrace}

\newif\if@nstr
\def\setstrfalse{\let\if@nstr=\iffalse}
\def\setstrtrue{\let\if@nstr=\iftrue}
\def\@nstr #1#2{
\def\@@nstr ##1#1##2##3\@@nstr{\ifx
\@nstr ##2\setstrfalse \else \setstrtrue \fi }
\@@nstr #2#1\@nstr \@@nstr}
\def\@separate#1|#2@{\setFront{#1}\setBack{#2}}
\def\lb#1\rb{\@nstr|{#1} \if@nstr \@separate#1 @ \twoSet{\@setFront}{\@setBack}%
\else \@separate |{#1 }@ \oneSet{\@setBack}\fi%
}
\def\setFront#1{\def\@setFront{#1}}
\def\setBack#1{\def\@setBack{#1}}
\def\Set#1{\lb{#1}\rb}
\def\oneBrk#1{\left\langle#1\right\rangle}
\def\twoBrk#1#2{\left\langle%
\vphantom{#2}#1\thinspace\right|\nolinebreak[3]\left.%
  #2%
  \vphantom{#1}%
  \right\rangle%
}
\def\brk<#1>{\@nstr|{#1} \if@nstr \@separate#1 @ \twoBrk{\@setFront}{\@setBack}%
\else \@separate |{#1 }@ \oneBrk{\@setBack}\fi%
}
\def\thmref#1{\normalfont{theorem}~\ref{#1}}
\def\lemref#1{\normalfont{lemma}~\ref{#1}}

\theoremstyle{plain}
\newtheorem{thm}{Theorem}[section]
\newtheorem{lem}[thm]{Lemma}
\newtheorem{cor}[thm]{Corollary}
\newtheorem{prop}[thm]{Proposition}
\newtheorem{defn}[thm]{Definition}

\theoremstyle{remark}
{}
{}
{}
{\newtheorem{rem}{Remark}[section]}

\@ifpackageloaded{bbm}{\newcommand{\Z}{{\mathbbm{Z}}}

}{
\newcommand{\Z}{{\mathbb{Z}}}

}
\makeatother

\begin{document}
\begin{abstract}
	In $\mathcal L$, the semilattice of faces of an $n$-cube, we count the 
	number of automorphisms of $\mathcal L$ that fix a given 
	subalgebra -- either pointwise or as a subalgebra. 
	By using M\"obius inversion we get a formula for the number of 
	derangements on the $n$-cube in terms of the M\"obius 
	function on the lattice of MR-subalgebras. We compute this 
	M\"obius function.
\end{abstract}
\maketitle
\section{Introduction}
We are interested in derangements of the $n$-cube,  ie the 
automorphisms that fix only the codimension zero face of the cube. Our approach is to 
consider the face-semilattice of the $n$-cube -- $\mathcal L_{n}$ -- and use the fact that 
this semilattice is a Metropolis-Rota implication algebra (MR-algebra), and the 
automorphism group of this algebra has a well-known structure. 
Since any automorphism that fixes a set $X$ of edges of the 
$n$-cube must also fix the entire MR-subalgebra of $\mathcal L_{n}$
generated by $X$, we restrict our study to 
MR-subalgebras of $\mathcal L$. 

We begin by fixing a subalgebra $A$ and consider the ways that $A$ 
might be fixed. 
\begin{defn}\label{def:fix}
	Let $A$ be a sub-MR-algebra of $\mathcal L$. Let $\phi$ be an 
	automorphism of $\mathcal L$. 
	\begin{enumerate}[(a)]
		\item $\phi$ \emph{freezes} $A$ iff $\phi\restrict 
		A=\text{id}\restrict A$. 
	
		\item $\phi$ \emph{fixes} $A$ iff $\phi[A]= A$. 
	\end{enumerate}
\end{defn}

This leads us to the following definition of two subgroups associated 
with $A$.

\begin{defn}\label{def:fixGrp}
	Let $A$ be a sub-MR-algebra of $\mathcal L$. 
	\begin{enumerate}[(a)]
		\item $\Fr(A)=\Set{\phi | \phi\text{ freezes }A}$; 
	
		\item $\Fix(A)=\Set{\phi | \phi\text{ fixes }A}$. 
	\end{enumerate}
\end{defn}

It is easy to see that both $\Fr(A)$ and $\Fix(A)$ are subgroups of $\Aut(\mathcal L)$. 
We want to count the size of each of these groups. We also want to 
determine the number of automorphisms that freeze $A$ only. 
This we will do by M\"obius inversion, as if 
\begin{align*}
	f(A) & =\card{\Fr(A)}  \\
	s(A) & =\card{\Set{\phi\in\Aut(\mathcal L) | \phi\restrict 
	A=\id_{A}\text{ and }\forall B>A\ \phi\restrict 
	B\not=\id_{B}}}  \\
	g(A) & =\card{\Fix(A)}  \\
	\intertext{then we see that}
	f(A) & =\sum_{A\subseteq B}s(B)  \\
	\intertext{so that, by M\"obius inversion,  we have }
	s(A) & =\sum_{A\subseteq B}\mu(A, B)f(B)
\end{align*}
where $\mu$ is the M\"obius function on the partial order of 
MR-subalgebras of $\mathcal L$. 
As a special case we get $s(\Set{\one})$ is the number of 
derangements on $\mathcal L$. 

So we will compute the functions $f$ and $g$ and the M\"obius function on the 
partial order of MR-subalgebras of $\mathcal L$. 
To evaluate $f$ and $g$ we compute orbits.

There is a study of derangements of $n$-cubes by Chen \& Stanley in 
\cite{CS:der}. That study concentrates on properties of signed 
permutations (which are the automorphisms of $\mathcal L$)
and obtains an alternative counting of the derangments. 

We begin by giving some basic background on cubic and MR implication 
algebras. The reader is referred to \cite{BO:eq} or 
\cite{MR:cubes} for a more thorough introduction to cubic and MR-algebras.

\subsection{Background Material}

\subsubsection{Cubic and MR-algebras}
We will give a brief introduction to the basic properties of cubic and 
MR-algebras. 
A \emph{cubic algebra} is an upper semilattice $\mathcal L$ with a binary operator $\Delta$ on $\mathcal 
L$ -- $\Delta(x, y)$ corresponding to reflection of $y$ through the 
centre of $x$ -- satisfying the following axioms:
\begin{enumerate}[a.]
    \item  if $x\le y$ then $\Delta(y, x)\join x = y$;
    
    \item  if $x\le y\le z$ then $\Delta(z, \Delta(y, x))=\Delta(\Delta(z, 
    y), \Delta(z, x))$;
    
    \item  if $x\le y$ then $\Delta(y, \Delta(y, x))=x$;
    
    \item  if $x\le y\le z$ then $\Delta(z, x)\le \Delta(z, y)$;
    
    \item[] Let $xy=\Delta(1, \Delta(x\join y, y))\join y$ for any $x$, $y$ 
    in $\mathcal L$. Then:
    
    \item  $(xy)y=x\join y$;
    
    \item  $x(yz)=y(xz)$;
\end{enumerate}

In fact the face poset of an $n$-cube, $\mathcal L_{n}$, is also an MR-algebra. We recall the pertinent 
details. 

\begin{defn}
    An \emph{MR-algebra} is a cubic algebra satisfying the MR-axiom:\\
    if $a, b<x$ then 
    \begin{gather*}
        \Delta(x, a)\join b<x\text{ iff }a\meet b\text{ does not exist.}
    \end{gather*}
\end{defn}

\begin{defn}\label{def:caret}
    Let $\mathcal L$ be a cubic algebra. Then for any $x, y\in\mathcal L$ 
    we define the (partial) operation $\caret$ (\emph{caret}) by:
    $$
        x\caret y=x\meet\Delta(x\join y, y)
    $$
    whenever this meet exists. 
\end{defn}

\begin{lem}
    If $\mathcal L$ is a cubic algebra then 
    $\mathcal L$ is an MR-algebra iff the caret operation is total.
\end{lem}
\begin{proof}
    See \cite{BO:fil} lemma 10 and theorem 12. 
\end{proof}

\begin{defn}
    Let $\mathcal L$ be a cubic algebra and $a, b\in\mathcal L$. Then
    \begin{align*}
        a\preceq b &\text{ iff }\Delta(a\join b, a)\le b\\
        a\simeq b &\text{ iff }\Delta(a\join b, a)=b.
    \end{align*}
\end{defn}

\begin{lem}
    Let $\mathcal L$, $a$, $b$ be as in the definition. Then
    $$
    a\preceq b\text{ iff }b=(b\join a)\meet(b\join\Delta(1, a)).
    $$
\end{lem}
\begin{proof}
    See \cite{BO:eq} lemmas 2.7 and 2.12.
\end{proof}

\begin{lem}
	Let $\mathcal L$ be a cubic algebra and $a\in\mathcal L$. If $b, 
	c\geq a$ then
	$$
	b\preceq c\iff b\le c.
	$$
\end{lem}
\begin{proof}
	If $b\le\Delta(b\join c, c)$ then we have 
	$a\le c$ and $a\le b\le \Delta(b\join c, c)$ and so 
	$b\join c = a\join\Delta(b\join c, a)\le c\join c=c$. 
\end{proof}

There are a number of representations of MR-algebras. For finite ones 
the principal three are as the face lattice of an $n$-cube; as the 
poset by signed subsets of $\Set{1, \dots, n}$; and as the poset of 
closed intervals of $\wp(\Set{1, \dots, n})$. We will consider the 
latter two briefly. 
\subsection{Signed Sets}

\begin{defn}\label{def:ss}
	Let $X$ be a set. 
	\begin{enumerate}[(a)]
		\item A \emph{signed subset} of $X$ is a pair $\brk<A_{1}, A_{2}>$ 
		where $A_{i}\subseteq X$ and $A_{1}\cap A_{2}=\emptyset$. 
	
		\item $\rsf S(X)$ is the collection of all signed subsets of $X$ 
		ordered by reverse pointwise inclusion. 
	\end{enumerate}
\end{defn}

$\rsf S(X)$ is an MR-algebra and if $\mathcal L$ is any finite 
MR-algebra with $\text{CoAt}(\mathcal L)$ its set of coatoms and 
$C\subset\text{CoAt}(\mathcal L)$ is such that 
$C\cup\Set{\Delta(\one, c) | c\in C}=\text{CoAt}(\mathcal L)$ and 
$C\cup\Set{\Delta(\one, c) | c\in C}=\emptyset$ then there is a 
canonical isomorphism of $\mathcal L$ with $\rsf S(C)$ -- see 
\cite{MR:cubes} for more details. 

As one application of this construction we have a simple homomorphism 
extension result for finite MR-algebras. 

\begin{prop}\label{prop:isoTwo}
    Let $\mathcal L$ be a finite MR-algebra. Let $\phi\colon\text{CoAt}(\mathcal 
    L)\to\text{CoAt}(\mathcal L)$ be a $\Delta(\one, \bullet)$-preserving 
    bijection. Then there is a canonical extension of $\phi$ to an 
    automorphism of $\mathcal L$. 
\end{prop}
\begin{proof}
    Let $a$ be an atom of $\mathcal L$ and let $C$ be the coatoms over $a$. 
    Then 
    $C\cup\Delta(\one, C)=\text{CoAt}(\mathcal L)$ and 
    $C\cup\Delta(\one, C)=\emptyset$. Likewise $\phi[C]$ has the same 
    properties since $\phi$ preserves $\Delta(\one, \bullet)$. From 
    \cite{CO:impl} this implies $\bigwedge\phi[C]$ exists and is an atom $a'$
    of $\mathcal L$. 
    Now $[a, \one]$ and $[a', \one]$ are isomorphic as Boolean 
    algebras by an extension of $\phi\restrict C$ and so we get 
    an extension of this mapping to an automorphism of $\mathcal L$. 
    This extends $\phi$ also as it has to preserve $\Delta(\one, \bullet)$. 
    
    An alternative way to view this proof is via the isomorphism sequence 
    $$
    \mathcal L\rTo^{\sim}\rsf S(C)\rTo^{\rsf S(\phi)}\rsf 
    S(\phi[C])\rTo^{\sim}\mathcal L
    $$
    where 
    $\rsf S(\phi)(\brk<A, B>)=\brk<\phi[A], \phi[B]>$. 
\end{proof}

\subsection{Implication Algebras}
Let $\mathcal I$ be an implication algebra (ie an upwards closed subset 
of a Boolean algebra). We define
$$
\rsf I(\mathcal I)=\Set{\brk<a, b> | a, b\in\mathcal I, a\join
b=\one\text{ and }a\meet b\text{ exists in }\mathcal I}
$$
ordered by 
$$
\brk<a, b>\le\brk<c, d>\text{ iff }a\le c\text{ and }b\le d. 
$$
This is a partial order that is an upper semi-lattice with join 
defined by
$$
\brk<a, b>\join\brk<c, d>=\brk<a\join c, b\join d>
$$
and a maximum element $\one=\brk<1, 1>$. 

We can also define a $\Delta$ function by
$$
\text{if }\brk<c, d>\le\brk<a, b>\text{ then }
\Delta(\brk<a, b>, \brk<c, d>)=\brk<a\meet(b\to d), b\meet(a\to c)>. 
$$

More properties of this construction are described in \cite{BO:fil}. 

\subsection{The Problem}
In \cites{BO:cubAutFin, BO:fil, Joe:thesis} the automorphism group of $\mathcal 
L$ was investigated.  In this paper we wish to consider automorphisms 
that fix an MR-subalgebra of $\mathcal L$. 

As described in the introduction we have the two groups 
$\text{Fr}(A)$ and $\text{Fix}(A)$ ($=\stab(A)$, the stabilizer of $A$) and to find the size of each of 
these groups we consider the orbit of $A$

\section{Orbits}
We will consider the natural group action of $\Aut(\mathcal L)$ on 
subalgebras. First we 
want to examine an invariant (the type) of an MR-subalgebra.

\begin{defn}\label{def:coat}
    Let $A$ be an MR-subalgebra of $\mathcal L_{n}$. 
    
    $\text{CoAt}_{n}$ is the set of coatoms of $\mathcal L_{n}$, 
    $\text{CoAt}(A)$ is the set of coatoms of $A$. 
    
    For each $a\in\text{CoAt}(A)$ we let $\rsf 
    C_{a}=\Set{c\in\text{CoAt}_{n} | a\le c}$ and $\Gamma_{A}=\Set{\rsf 
    C_{a} | a\in \text{CoAt}(A)}$. 
\end{defn}

We notice that $a\in\text{CoAt}(A)$ implies $\Delta a\in\text{CoAt}(A)$ 
and $\rsf C_{\Delta a}=\Delta\rsf C_{a}$. Thus for each 
$1\le i\le n$ there are an even number of $\rsf C_{a}$'s of size $i$. 
Let $t_{i}$ be such that $2t_{i}=\card{\Set{a\in\text{CoAt}(A) | 
\card{\rsf C_{a}}=i}}$. 

\begin{defn}\label{def:typeA}
    The \emph{type} of an MR-subalgebra $A$ of $\mathcal L_{n}$ is the 
    sequence $\text{tp}(A)=\brk<t_{i} | 1\le i\le n>$. 
\end{defn}

The action we are considering is the evaluation action of
$\Aut(\mathcal L)$ on the partial order of 
MR-subalgebras of $\mathcal L_{n}$. 
Thus we know that 
\begin{align*}
	\card{\Fix(A)}&=\card{\Aut(\mathcal L)}/\card{\text{Orb}(A)}\\
	&= 2^{n}n!/\card{\text{Orb}(A)}
\end{align*}
and we therefore want to compute the size of the orbits.

\begin{lem}\label{lem:iOne}
	If $A$ and $B$ are two MR-subalgebras on $\mathcal L$ in the same 
	orbit, then $\text{tp}(A)=\text{tp}(B)$. 
\end{lem}
\begin{proof}
	This is clear. 
\end{proof}

\begin{lem}\label{lem:autExtend}
	Let $c_{1}$ and $c_{2}$ be two coatoms of $\mathcal L$ and 
	$\phi\colon]\leftarrow, c_{1}]\to ]\leftarrow, c_{2}]$ be an 
	isomorphism. Then $\phi$ extends to an automorphism of $\mathcal L$. 
\end{lem}
\begin{proof}
	The coatoms of $]\leftarrow, c_{1}]$ are of the form $c\meet c_{1}$ 
	where $c\not= c_{1}$ and $c\not=\Delta(\one, c_{1})$ -- see 
	\cite{MR:cubes}. Furthermore $c\not=d$ implies $c\meet 
	c_{1}\not=d\meet c_{1}$. So we define the extension by defining what 
	happens on the coatoms:
	\begin{align*}
		\phi'(c) & =c'\text{ where }\phi(c\meet c_{1})=c'\meet c_{2}  \\
		\phi'(c_{1}) & =c_{2}  \\
		\phi'(\Delta(\one, c_{1})) & =\Delta(\one, c_{2}). 
	\end{align*}
	From \cite{CO:impl} we know that $\Delta(c_{i}, c_{i}\meet 
	c)=c_{i}\meet\Delta(\one, c)$ and so $\phi'(\Delta(\one, c))=\Delta(\one, \phi'(c))$. 
	
	From this we extend $\phi'$ to an automorphism of $\mathcal L$ as usual. 
\end{proof}

\begin{lem}\label{lem:iTwo}
	Let $v_{1}$ and $v_{2}$ be two elements of $\mathcal L$ of the same 
	co-rank. Then there is an automorphism of $\mathcal L$ taking 
	$v_{1}$ to $v_{2}$. 
\end{lem}
\begin{proof}
	We proceed by induction on co-rank. Suppose that $v_{1}$ and $v_{2}$ 
	are coatoms. There are two cases. 
	\begin{enumerate}[{Case }1:]
		\item $v_{1}=\Delta(\one, v_{2})$ -- in this case just take 
		$\Delta(\one, \bullet)$ as the automorphism. 
	
		\item $v_{1}\not=\Delta(\one, v_{2})$ -- then $v_{1}\meet v_{2}$ 
		exists and so we have a Boolean algebra $[v, \one]$, where $v$ is any 
		vertex below $v_{1}\meet v_{2}$. As $v_{1}$ and $v_{2}$ are coatoms 
		of $[v, \one]$ there is a Boolean automorphism of $[v, \one]$ 
		taking $v_{1}$ to $v_{2}$. As usual this extends to an automorphism 
		of $\mathcal L$. 
	\end{enumerate}
	
	Now let $c_{i}\geq v_{i}$ be two coatoms and let $\phi$ be an 
	automorphism of $\mathcal L$ that takes $c_{1}$ to $c_{2}$. 
	In $]\leftarrow, c_{2}]$ the co-rank of $v_{2}$ equals the co-rank of 
	$\phi(v_{1})$ and is one less than the co-rank of $v_{i}$ in $\mathcal 
	L$. By induction, there is an automorphism $\psi$ of $]\leftarrow, c_{2}]$ 
	that takes $\phi(v_{1})$ to $v_{2}$. By \lemref{lem:autExtend} this 
	extends to an automorphism $\psi'$ of $\mathcal L$. Then we have 
	$\psi'\phi(v_{1})=v_{2}$ as desired. 
\end{proof}

\begin{thm}\label{thm:iTwo}
	Suppose that $A$ and $B$ are two MR-subalgebras on $\mathcal L$ with  
	$\text{tp}(A)=\text{tp}(B)$. Then there is an automorphism of $\mathcal L$ that 
	takes $A$ to $B$. 
\end{thm}
\begin{proof}
	Let $\text{tp}(A)=\text{tp}(B)=\brk<t_{i} | 1\le i\le n>$. 
	Let $r=\sum_{i=1}^{n}it_{i}$. This must be the corank of a vertex in $A$ 
	or $B$. Let $k=\sum_{i=1}^{n}t_{i}$ -- this is the dimension of $A$ (ie 
	$A\simeq \mathcal L_{k}$). 
	
	First we find $v_{1}\in A$ and $v_{2}\in B$ having co-rank $r$, and 
	find an automorphism of $\mathcal L$ that takes $v_{1}$ to $v_{2}$. 
	
	Now the $A$-atoms over $v_{1}$ go to an antichain in $[v_{2}, \one]$ 
	that induces a partition of the $\mathcal L$-covers of $v_{2}$. 
	Likewise the $B$-atoms of $[v_{2}, \one]$ 
	induce a partition of the $\mathcal L$-covers of $v_{2}$. Since 
	$\text{tp}(A)=\text{tp}(B)$ these partitions are similar and so there is a 
	permutation of the $\mathcal L$-covers of $v_{2}$ taking the first 
	antichain to the second. This induces a Boolean automorphism $\psi$ of $[v_{2}, 
	\one]$. 
	
	Now, let $a_{1}\le v_{1}$ be an $\mathcal L$-vertex and 
	$a_{2}=\phi(a_{1})\le v_{2}$. Then the automorphism $\psi$ can be 
	extended to a Boolean automorphism $\psi'$ of $[a_{2}, \one]$ and to 
	an automorphism $\widehat{\psi'}$ of $\mathcal L$. Then we have 
	$\widehat{\psi'}\phi[A]=B$. 
\end{proof}

From the lemma we need only count the number of 
ways we get the same type. Suppose that $\text{tp}(B)=\text{tp}(A)$ and $d$ is a 
vertex of $B$.  
We can think of a set of atoms of a subalgebra of $[d, \one]$ as a 
partition of the $\mathcal L$ covers of $d$. If $\cork(d)=r$ then 
$\rk[d, \one]=r$ so these correspond to partitions of $r$. 
$\At(A)$ gives one such partition $\Pi$ and we need to count the number of 
similar partitions. 

We also need to recall that each MR-subalgebra of rank $k$ has $2^{k}$ 
vertices and so 
the size of the orbit is 
$$
	\card{\Set{ d\in\mathcal L | \cork(d)=r}}/2^{k}\times\text{the 
	number of partitions of $r$ similar to }\Pi. 
$$

The number $\card{\Set{ d\in\mathcal L | \cork(d)=r}}$ is 
known to be 
$$
\card{\Set{ d\in\mathcal L | 
\cork(d)=r}}=2^{r}\binom{n}{r}
$$
-- see \cite{Benn} for example. The number of partitions of $r$ similar 
to $\Pi$ is well known to be 

$$
\frac{r!}{\prod_{i=1}^{n}(i!)^{t_{i}}t_{i}!}. 
$$
Thus 
the size of the orbit is 
$$
2^{r-k}\binom{n}r\frac{r!}{\prod_{i=1}^{n}(i!)^{t_{i}}t_{i}!}. 
$$

From this we infer that 
$$
\card{\Fix(A)}=\frac{2^{n}n!\prod_{i=1}^{n}(i!)^{t_{i}}t_{i}!}{2^{r-k}\frac{n!}{(n-r)!}}=
2^{n+k-r}(n-r)!\prod_{i=1}^{n}(i!)^{t_{i}}t_{i}!. 
$$

There are two ways to compute the size of $\Fr(A)$ -- one directly 
and another
by noticing that we have a group homomorphism 
$$
\rho\colon \varphi\mapsto \varphi\restrict A
$$
from $\Fix(A)$ to $\Aut(A)$ with kernel $\Fr(A)$, and so we can 
compute the image of this homomorphism. 

We will do both as each gives a viewpoint on automorphisms that we find 
interesting. 

We can partition $\text{CoAt}(A)$ into the sets $\Gamma_{i}=\Set{a | 
\cork(a)=i}$. We notice that $\card{\Gamma_{i}}=2t_{i}$ and 
$a\in\Gamma_{i}$ implies $\Delta a\in\Gamma_{i}$.

\begin{lem}\label{lem:extend}
    Let $\psi\in\Aut(A)$. Then $\psi\in\text{Im}(\rho)$ iff 
    $\psi[\Gamma_{i}]=\Gamma_{i}$ for all $i$. 
\end{lem}
\begin{proof}
    Let $\psi=\rho(\Psi)$ for some $\Psi\in\text{Stab}(A)$. The
    $\Psi$ preserves corank and takes $A$ to $A$,  so it must take 
    $\Gamma_{i}$ into itself. As it is one-one on $\Gamma_{i}$ it is 
    also onto. 
    
    For the converse,  suppose that $\psi\in\Aut(A)$ is such that 
    $\psi[\Gamma_{i}]=\Gamma_{i}$ for all $i$.
    
    Let $N_{i}\subseteq\Gamma_{i}$ be maximal $\Delta$-independent. 
    
    We claim that $\rsf C_{N_{i}}$ (ie the set of coatoms above some 
    element of $N_{i}$) is also $\Delta$-independent in $\rsf 
    C_{\Gamma_{i}}$. Indeed,  if $c\in\rsf C_{N_{i}}$ is greater than 
    $a\in N_{i}$ and $\Delta c\geq b\in N_{i}$ then $c\in\rsf 
    C_{a}\cap\rsf C_{\Delta b}$. As $\Delta b\in N_{i}$ and 
    $a\not=\Delta b$ implies $\rsf C_{a}\cap\rsf C_{\Delta 
    b}=\emptyset$ we must have $a=\Delta b$ -- contradicting the 
    assumption that $N_{i}$ is $\Delta$-independent. 
    
    Let $N_{i}'=\psi[N_{i}]\subseteq \Gamma_{i}$. This is also 
    $\Delta$-independent and so is $\rsf C_{N_{i}'}$. 
    
    $\rsf C_{N_{i}}$ is the disjoint union of the set $\Set{\rsf 
    C_{a} | a\in N_{i}}$ and we know that $\card{\rsf 
    C_{a}}=\card{\rsf C_{\psi(a)}}= i$ so we can find a bijection
    from $\rsf C_{N_{i}}$ to $\rsf C_{N_{i}'}$ that takes 
    $\rsf C_{a}$ to $\rsf C_{\psi(a)}$. 
    
    Finally we can patch these bijections to get a bijection
    between the disjoint union of the set $\Set{\rsf C_{N_{i}}| 1\le 
    i\le n}$ and the disjoint union of the set $\Set{\rsf C_{N_{i}'}| 1\le 
    i\le n}$. This is now a bijection between two maximal 
    $\Delta$-indpendent subsets of $\text{CoAt}(\mathcal L)$ and so 
    extends to an automorphism of $\mathcal L$. 
    
    We note that it takes $a\in\Gamma_{i}$ to $\psi(a)$ and so it 
    restricts to $\psi$. 
\end{proof}

\begin{lem}\label{lem:part}
    Let $P$ be a partition of $\text{CoAt}_{n}$ such that for all 
    $X\in P$\ $\Delta[X]=X$. Let 
    $\Aut_{P}=\Set{\varphi\in\Aut(\mathcal L_{n})| \forall X\in P\ 
    \varphi[X]=X}$. Then 
    $$
    \Aut_{P}\simeq\prod_{X\in P}\Aut(\mathcal L_{\card X/2}). 
    $$
\end{lem}
\begin{proof}
We first observe that if $\Delta[X]=X$ for some set of coatoms $X$ of 
$\mathcal L_{n}$ then 
$A_{X}=\Set{x\in\mathcal L_{n} | \text{the coatoms above }x\text{ are all in 
}X}$ is an MR-subalgebra isomorphic to $\mathcal L_{\card X/2}$ -- by 
taking $X'\subseteq X$ so that $X'\cap\Delta[X']=\emptyset$ and 
$X'\cup\Delta[X']=X$ and the mapping 
$x\mapsto \brk<\Set{z\in X'| x\le z}, \Set{z\in\Delta[X'] | x\le z}>$ 
is a cubic isomorphism to the algebra of signed subsets of $X'$, ie 
$\mathcal L_{\card X/2}$. 

Thus, if $\varphi\in\Aut_{P}$ then $\varphi\restrict A_{X}$ is a cubic 
automorphism of $A_{X}$.  

This mapping is onto -- as if $\brk<\varphi_{X} | X\in P>\in \prod_{X\in 
P}\Aut(A_{X})$, then we consider the mapping defined on 
$\text{CoAt}_{n}$ by 
$$
x\mapsto \varphi_{X}(x) \text{ if }x\in X. 
$$
This is a $\Delta$-preserving mapping of the coatoms and so lifts to an 
automorphism of $\mathcal L_{n}$, which restricts to $\varphi_{X}$ on 
each $A_{X}$. 
\end{proof}

From these two lemmas we see that if $P=\Set{\Gamma_{i} | 1\le i\le n}$
is the partition of $\text{CoAt}(A)$ described above, then 
$$
\card{\text{Im}(\rho)}=\card{\Aut_{P}}=\prod_{i}\card{\Aut(\mathcal 
L_{t_{i}})}= \prod_{i}2^{t_{i}}t_{i}!= 2^{k}\prod_{i}t_{i}!.  
$$

Thus we have 
$$
\card{\Fr(A)}=\frac{2^{n+k-r}(n-r)!\prod_{i=1}^{n}(i!)^{t_{i}}t_{i}!}{2^{k}\prod_{i}t_{i}!}=
2^{n-r}(n-r)!\prod_{i=1}^{n}(i!)^{t_{i}}. 
$$

Here is another way to see this result -- consider the sets $\rsf C_{a}$ 
as $a$ varies over $\text{CoAt}(A)$, and 
$D=\text{CoAt}_{n}\setminus\bigcup_{a}\rsf C_{a}$. Notice that 
$\Delta[D]=D$, and $\Delta[\rsf C_{a}]=\rsf C_{\Delta a}$ for all $a$.  

Then $\varphi$ freezes $A$ iff $\varphi[\rsf C_{a}]=\rsf C_{a}$ for all 
$a\in\text{CoAt}(A)$. This of course, implies $\varphi\restrict D$ is 
a $\Delta$-preserving mapping from $D$ to itself. 

Also we must have $\varphi\restrict\rsf C_{\Delta 
a}=\Delta\bigl(\varphi\restrict\rsf C_{a}\bigr) \Delta$ -- as $\varphi$ 
preserves $\Delta$. 

Hence, if $M$ is a maximal $\Delta$-independent set of coatoms of $A$ 
then $\varphi$ is completely determined by its action on $\rsf C_{a}$ 
for $a\in M$ and its action on $D$. $\varphi$ can be any permutation of 
$\rsf C_{a}$ (for $a\in M$) and (as above) any $\Delta$-preserving 
bijection of $D$ -- of which there are $2^{n-r}(n-r)!$ such mappings. Hence there are
$$
2^{n-r}(n-r)!\prod_{i}(i!)^{t_{i}}
$$
such $\varphi$ -- as computed above.

\section{The M\"obius function on Implication lattices}
As the first step in computing the M\"obius function on the poset of 
MR-subalgebras we will look at implication sublattices of a Boolean 
algebra. The MR-subalgebra computation will then be reduced to this case.

Let $B$ be a finite Boolean algebra. We will assume that $B\simeq\wp(n)$. 
\begin{defn}\label{def:implsub}
	An \emph{implication subalgebra} of $B$ is a subset closed under $\to$. 
	
	An \emph{implication sublattice} of $B$ is a  subset closed under 
	$\to$  and $\wedge$. 
\end{defn}

In \cite{BO:implMob} we found that this function is given by the 
formula
\begin{align}
\mu(A, B)&=
(-1)^{n-k}(n-k)!\prod_{i=1}^{n}\left[(-1)^{i-1}(i-1)!\right]^{t_{i}}\\
&= 
(-1)^{n-d}(n-k)!\prod_{i=1}^{n}\left[(i-1)!\right]^{t_{i}}
\end{align}
as $\sum_{i=1}^{n}(i-1)t_{i}=k-d$,  and 
\begin{equation}\label{eq:muOne}
\mu(\Set{\one}, B)=n!. 
\end{equation}

\section{The M\"obius function on MR-subalgebras}
The way we will compute the M\"obius function on the poset of 
MR-subalgebras is similar to that of the last section. 

We begin by showing how $\mu(M,N)$ can always be determined by knowing 
$\mu(\Set{\one}, \mathcal L_{k})$ for all $k$. 

Then we 
represent the poset in terms of implication sublattices of 
$[\zero, \one]$ together with some extra information. Then we define 
a closure operator on this new representation and finally reduce the 
problem to the poset of implication sublattices of $[\zero, \one]$.

Let $A$ and $B$ be two MR-subalgebras of $\mathcal L$. As $B$ is 
already an MR-algebra and $B$ so is the face lattice of a (possibly) smaller 
cube, we may assume that $B=\mathcal L$. We consider a reduction 
showing that we may also assume that $A=\Set{\one}$. 

Fix an atom $\zero$ of $\mathcal L$ below some $A$-atom $a_{A}$. 
Let $A\subseteq C\subseteq \mathcal L$ be any intermediate subalgebra. 
Then $C$ is determined by knowing $C\cap[a_{A}, 1]$ and 
$C\cap]\leftarrow, a_{A}]$. This shows us that 
\begin{multline*}
[A, \mathcal L]\simeq\Set{ B | A\cap [a_{A}, \one]\subseteq B\subseteq 
[a_{A}, \one]\text{ is a Boolean algebra}}\\
\times\Set{ C | C\text{ is 
an MR-subalgebra of }]\leftarrow, a_{A}]}. 
\end{multline*}
The first factor is well-known as a poset of partitions. The second 
factor is our poset with $A=\Set{\one}$. 

So we will assume that $A=\Set{\one}$. 

\subsection{Locator Pairs}
We consider MR-subalgebras in $[\Set{\one}, \mathcal L]$ using 
implication sublattices of $[\zero, \one]$. In this section we develop 
a way of describing subalgebras that leads to a clearer picture of the 
partial order. 

\begin{defn}\label{def:locator}
	A \emph{locator-pair} is a pair of $\brk<c, B>$ where
	$B$ is an implication sublattice of $[\zero, \one]$ and 
	$c\geq \min B$. 
\end{defn}

Locator pairs will be used to facilitate counting. 

\begin{lem}\label{lem:locAlg}
	Locator pairs correspond to subalgebras of $\mathcal L$. 
\end{lem}
\begin{proof}
	Let $A$ be any subalgebra, let $a_{A}$ be a vertex of $A$. 
	Let $A_{*}=\beta_{\zero}[A]$ and $c=\zero\join a_{A}$. Then 
	$\brk<c, A_{*}>$ is a locator-pair. We can recover $A$ from this 
	locator-pair by noting that if 
	$a_{A}=\min A_{*}$ then $a=\Delta(c, a_{A})$ is a vertex of $A$ and we 
	can move $A_{*}$ to a g-filter of $A$ using the mapping 
	$$
	x\mapsto (x\join a)\meet(\Delta(\one, x)\join a)
	$$
	(as usual), and as this is a g-filter for $A$ we can recover $A$. 
\end{proof}

The implication sublattice $A_{*}$ is uniquely determined by $A$, but 
the other element is not -- as we can choose many $A$-vertices. 

\begin{defn}\label{def:locate}
	Let $A$ be an MR-subalgebra of $\mathcal L$. A locator-pair $\brk<c, 
	B>$ that determines $A$ is said to \emph{locate} $A$. 
\end{defn}
By an abuse of notation we will often write this as $\brk<c, B>=A$. 
 
\begin{defn}\label{def:simeq}
	Let $\brk<c_{1}, B_{1}>$ and $\brk<c_{2}, B_{2}>$ be two locator-pairs. 
	Let $A_{i}$ be the corresponding MR-subalgebra of $\mathcal L$. Then
	\begin{align*}
		\brk<c_{1}, B_{1}>\le \brk<c_{2}, B_{2}> & \text{ iff }A_{1}\subseteq A_{2}  \\
		\brk<c_{1}, B_{1}>\simeq\brk<c_{2}, B_{2}> & \text{ iff }A_{1}= A_{2}. 
	\end{align*}
\end{defn}

It is easy to see that $\simeq$ is an equivalence relation. We want 
to characterize $\le$ on locator-pairs more carefully. 

\begin{defn}\label{def:plusSuba}
	Let $c, d$ be in $\mathcal L$. Then 
	$c+_{a} d$ is the Boolean sum of $c\join a$ and $d\join a$ in $[a, \one]$. 
\end{defn}

Note that in a Boolean algebra we have $c+_{a}d=(c+d)\join a$. Hence 
in a cubic algebra or implication algebra, if $a_{1}\le a_{2}$ then 
$c+_{a_{2}}d=(c+_{a_{1}}d)\join a_{2}$. 

We need the following technical lemma 
\begin{lem}\label{lem:plusLem}
	Let $c_{1}$ and $c_{2}$ be two elements of $\mathcal L$ such that 
	$c_{1}\meet c_{2}$ exists. Let $a\le c_{1}\meet c_{2}$. Then
	\begin{enumerate}[(a)]
		\item  if $a\le b\le c_{1}$ then 
		$$
		\Delta(c_{1}, b)\join\Delta(c_{2}, a)=\Delta(c_{1}, 
		b)\join\Delta(c_{2}\join b, b);  
		$$
	
		\item 
		if $a\le x$ then 
		$$
		\Delta(x\join c_{1}\join c_{2}, \Delta(x\join c_{1}, x)\join 
		\Delta(x\join c_{2}, x))=c_{1}+_{x}c_{2}. 
		$$
	\end{enumerate}
\end{lem}
\begin{proof}
	Without loss of generality we may work in an interval algebra and 
	take $a=[0, 0]$. Then we have
	$c_{i}=[0, c_{i}]$, $b=[0, b]$ and $x=[0, x]$. 
	\begin{enumerate}[(a)]
		\item 
		\begin{align*}
			\Delta(c_{1}, b)\join\Delta(c_{2}, a) & =\Delta([0, c_{1}], [0, 
			b])\join\Delta([0, c_{2}], [0, 0])  \\
			 & =[c_{1}\meet\comp b, c_{1}]\join[c_{2}, c_{2}]  \\
			 & =[c_{1}\meet c_{2}\meet\comp b, c_{1}\join c_{2}].   \\
			\Delta(c_{1}, b)\join\Delta(c_{2}\join b, b) & =
			\Delta([0, c_{1}], [0, b])\join\Delta([0, c_{2}\join b], [0, b])\\
			 & =[c_{1}\meet\comp b, c_{1}]\join[c_{2}\meet\comp b, c_{2}\join b]  \\
			 & =[c_{1}\meet c_{2}\meet\comp b, c_{1}\join c_{2}\join b]  \\
			 & =[c_{1}\meet c_{2}\meet\comp b, c_{1}\join c_{2}] && \text{ as 
			 }c_{1}\geq b. 
		\end{align*}
	
		\item 
		\begin{align*}
			\Delta(x\join c_{1}, x) & =\Delta([0, x\join c_{1}], [0, x])  \\
			& =[c_{1}\meet\comp x, c_{1}\join x]  \\
			\Delta(x\join c_{1}\join c_{2}, \Delta(x\join c_{1}, x)\join & \\ 
		\Delta(x\join c_{2}, x)) & =\Delta([0, x\join c_{1}\join c_{2}], 
		[c_{1}\meet\comp x, c_{1}\join x]\join[c_{2}\meet\comp x, c_{2}\join x])  \\
			 & = \Delta([0, x\join c_{1}\join c_{2}], 
		[c_{1}\meet c_{2}\meet\comp x, c_{1}\join c_{2}\join x]) \\
			 & =[0, (x\join c_{1}\join c_{2})\meet(x\join \comp{c_{1}}\join 
			 \comp{c_{2}})]  \\
			 & =[0, x\join(c_{1}+c_{2})] \\
			 & =[0, x]\join[0, c_{1}+c_{2}]  \\
			c_{1}+_{x}c_{2} & =x\join(c_{1}+_{a}c_{2}) \\
			c_{1}+_{a}c_{2} & =(c_{1}\join c_{2})\meet(c_{1}a\join c_{2}a)  \\
			 & =[0, c_{1}\join c_{2}]\meet[0, \comp{c_{1}}\join\comp{c_{2}}]\\
			 & = [0, c_{1}+c_{2}]. 
		\end{align*}
	\end{enumerate}
\end{proof}

\begin{thm}\label{thm:plusThm}
	Let $\brk<c_{1}, B_{1}>$ locate $A_{1}$ and $\brk<c_{2}, B_{2}>$ 
	locate $A_{2}$. Then 
	$$
		A_{1}\subseteq A_{2}\text{ iff }B_{1}\subseteq B_{2}\text{ and 
		}c_{1}+_{a_{1}}c_{2}\in B_{2}
	$$
	where $a_{i}=\min B_{i}$. 
\end{thm}
\begin{proof}
	Suppose that $A_{1}\subseteq A_{2}$. Then $B_{i}$ is obtained from 
	$A_{i}$ as the image of the mapping $x\mapsto\Delta(x\join \zero, x)$ 
	and so clearly $B_{1}\subseteq B_{2}$. 
	
	The locator $c_{i}$ has the property that $\Delta(c_{i}, a_{i})$ is 
	an atom of $A_{i}$. As $A_{1}\subseteq A_{2}$ this implies 
	$\Delta(c_{1}, a_{1})\join\Delta(c_{2}, a_{2})\in A_{2}$. Hence 
	$\Delta(\Delta(c_{1}, a_{1})\join\Delta(c_{2}, a_{2})\join a_{2}, 
	\Delta(c_{1}, a_{1})\join\Delta(c_{2}, a_{2}))\in B_{2}$.
	Now
	\begin{align*}
		\Delta(c_{1}, a_{1})\join\Delta(c_{2}, a_{2})\join a_{2} & =
		\Delta(c_{1}, a_{1})\join c_{2} \\
		&\qquad\qquad\qquad\qquad \text{as }\Delta(c_{2}, a_{2})\join 
		a_{2}=c_{2}\\
		 & \geq  \Delta(c_{1}, a_{2})\join c_{2}  \\
		 & = c_{1}\join c_{2}\\
		\text{As }\qquad\qquad \Delta(c_{1}, a_{1})\join c_{2} & \le c_{1}\join 
		c_{2}   \\
		\intertext{we have } 
		\Delta(c_{1}, a_{1})\join\Delta(c_{2}, a_{2})\join a_{2} & 
		=c_{1}\join c_{2}. \\
		\Delta(c_{1}, a_{1})\join\Delta(c_{2}, a_{2}) & = 
		\Delta(c_{1}, a_{1})\join\Delta(c_{2}\join a_{1}, a_{1}) \\
		& \qquad\qquad\qquad\qquad\text{by 
		\lemref{lem:plusLem} (a). } \\
		\intertext{Therefore }
		\Delta(\Delta(c_{1}, a_{1})\join\Delta(c_{2}, a_{2})\join a_{2}, & \\
	\Delta(c_{1}, a_{1})\join\Delta(c_{2}, a_{2})) &= 
	\Delta(c_{1}\join c_{2}\join a_{1}, 
	\Delta(c_{1}, a_{1})\join\Delta(c_{2}\join a_{1}, a_{1})) \\
	&= c_{1}+_{a_{1}}c_{2} \\
	& \qquad\qquad\qquad\qquad\text{ by \lemref{lem:plusLem} (b). }
	\end{align*}
	
	Now let us suppose that $B_{1}\subseteq B_{2}$ and 
	$c_{1}+_{a_{1}}c_{2}\in B_{2}$. It suffices to show that if $x\in 
	B_{1}$ then $\Delta(x\join\Delta(c_{1}, a_{1}), x)\in A_{2}$ -- as 
	the set of such elements forms a set that generates $A_{1}$.
	
	First we note that $x\geq a_{1}$ so that $x\join\Delta(c_{1}, 
	a_{1})=x\join c_{1}$. Also we know
	$x\in B_{2}$ so that $\Delta(x\join c_{2}, x)= 
	\Delta(x\join\Delta(c_{2}, a_{2}), x) \in A_{2}$ and furthermore 
	$\Delta(\Delta(x\join c_{1}, x)\join\Delta(x\join c_{2}, x), 
	\Delta(x\join c_{2}, x)) = \Delta(x\join c_{1}, x)$.  Thus it suffices to 
	show that $\Delta(x\join c_{1}, x)\join\Delta(x\join c_{2}, x)\in 
	A_{2}$.  
	
	For this it is sufficient to show that the preimage over $a_{2}$ is in 
	$B_{2}$ -- \\
	i.e. $\Delta(\Delta(x\join c_{1}, x)\join\Delta(x\join c_{2}, 
	x)\join a_{2}, \Delta(x\join c_{1}, x)\join\Delta(x\join c_{2}, x)) \in B_{2}$. 
	\begin{align*}
		\Delta(x\join c_{1}, x)\join\Delta(x\join c_{2}, 
	x)\join a_{2} & = x\join c_{1}\join c_{2}  \\
		\intertext{so that }
		\Delta(\Delta(x\join c_{1}, x)\join\Delta(x\join c_{2}, 
	x)\join a_{2}, & \\
	\Delta(x\join c_{1}, x)\join\Delta(x\join c_{2}, 
	x)) & =\Delta(x\join c_{1}\join c_{2}, \Delta(x\join c_{1}, x)\join\Delta(x\join c_{2}, x))  \\
		 & =c_{1}+_{x}c_{2} \qquad\qquad \text{by \lemref{lem:plusLem}}  \\
		 & =(c_{1}+_{a_{1}}c_{2})\join x \\
		 & \in B_{2} \qquad\qquad\text{ as }c_{1}+_{a_{1}}c_{2}\in B_{2}\text{ and 
		 }x\in B_{1}\subseteq B_{2}. 
	\end{align*}
\end{proof}

\subsection{A closure operator}
We will finally compute the M\"obius function we want through an 
appeal to the following theorem about closure operators  -- see 
\cite{INC} Proposition 2.1.19. 
\begin{thm}\label{thm:closure}
	Let $X$ be a locally finite partial order and $x\mapsto \comp x$ be 
	a closure operator on $X$. Let $\comp X$ be the suborder of all 
	closed elements of $X$ and $y$, $z$ be in $X$. Then 
	$$
		\sum_{\comp x=\comp z}\mu(y, x)=
		\begin{cases}
			\mu_{\comp X}(y, \comp z) & \text{ if }y\in\comp X  \\
			0 & \text{ otherwise. }
		\end{cases}
	$$
\end{thm}
\begin{proof}
	 See \cite{INC}. 
\end{proof}

There is several closure operators of interest that naturally apply to 
locator-pairs by modifying the second component. We will consider 
only one of them.

If $\brk<c, B>$ is a 
locator-pair we define $B^{*}$ to be the subalgebra of $[\min B, 1]$ 
generated by $\Set c\cup B$. Then we have 
$\comp{\brk<c, B>}=\brk<c, B^{*}>$. 

\begin{lem}\label{lem:closureUp}
	$\brk<c, B>\mapsto\comp{\brk<c, B>}$ is a closure operator on 
	locator-pairs. 
\end{lem}
\begin{proof}
	Since $B^{**}=B^{*}$ we trivially have 
	$\comp{\comp{\brk<c, B>}}=\comp{\brk<c, B>}$. 
	
	${\brk<c, B>}\le \comp{\brk<c, B>}$ as $B\subseteq B^{*}$ and 
	$c+_{b}c=b\in B^{*}$. 
	
	If $\brk<c_{1}, B_{1}>\le \brk<c_{2}, B_{2}>$ then 
	$B_{1}\subseteq B_{2}$ and so $B_{1}\subseteq B_{2}^{*}$. 
	Also $c_{1}+_{b_{1}}c_{2}\in B_{2}^{*}$ and therefore 
	$c_{1}=(c_{1}+_{b_{1}}c_{2})+_{b_{1}}c_{2}$ is in $B_{2}^{*}$. Hence
	$B_{1}^{*}\subseteq B_{2}^{*}$. ($b_{i}=\min B_{i}$). 
\end{proof}

\begin{lem}\label{lem:clos22}
    $\brk<c, B>$ is closed iff $\brk<c, B>\simeq\brk<\min B, B>$. 
\end{lem}
\begin{proof}
    Let $b=\min B$. 
    
It is clear that $\brk<c, B>$ is closed iff $c\in B$. Also, 
if $\brk<c, B>\simeq \brk<c', B>$ and $c'\in B$ then $c\in B$ -- since 
$c=(c+_{b}c')+_{b}c'$. Therefore $\brk<c, B>\simeq \brk<b, B>$ implies 
$c\in B$ and so $\brk<c, B>$ is closed.

Conversely, if $c\in B$ then $c+_{b}b= c\in B$ and so $\brk<c, 
B>\simeq \brk<b, B>$. 
\end{proof}

This lemma tells us that the poset of closed pairs is the same as the 
poset of implication sublattices of $[0, 1]$ -- that we have considered 
elsewhere \cite{BO:ImpMob}. 

We also note that $\brk<1, \Set1>$ is closed and locates $\Set1$, and 
$\brk<0, [0, 1]>$ is closed and locates $\mathcal L_{n}$.

\section{Getting the M\"obius Function}
We need to count the subalgebras whose closure is $\mathcal L$ rather 
carefully. 

We see that $\comp{\brk<c, B>}=\mathcal L$ iff $\min B=0$ and 
$B^{*}=[0,1]$. Thus $B$ is actually a Boolean subalgebra of $[0,1]$. 
Let $a_{1},\dots, a_{m}$ be the atoms of $B$.  The the atoms of 
$B^{*}$ are the non-zero elements of 
$$
\Set{a_{i}\meet c | 1\le i\le m}\cup\Set{a_{i}\meet\comp c|1\le i\le m}. 
$$
As these must be the atoms of $[0,1]$ we see that every atom of $B$ 
must be either a $[0,1]$-atom or the join of two such atoms. 

Let $k$ be the number of $B$-atoms that are also $[0,1]$-atoms and 
$\ell=m-k$.  Then we have $k+2\ell=n$ and the pair $\brk<c,B>$ is 
determined by the arrangement of $[0,1]$-atoms $\brk<S, P>$ where
\begin{align*}
S &= \Set{ a | a\text{ is a }B\text{-atom}}\\
P &= \Set{{\Set{a,b}} | a\neq b\text{ and }a\join b\text{ is a }B\text{-atom}}. 
\end{align*}

The number $k+l$ is the dimension of $B$ and the next lemma shows 
that this naturally determines a partition of the locators we are 
interested in. We also note that $k+l=m$ and $k+2\ell=n$ iff
$k=2m-n$ and $\ell=n-m$, so the pair is determined by the dimension of $B$ 
-- as $n$ is fixed. As we need $k,l\geq 0$ this also implies 
$n\geq m\geq\lceil\frac{n}2\rceil$. 

\begin{lem}\label{lem:atomsAreAll}
	The intervals $\mathbbm I_{1}=[\brk<\one, \Set{\one}>, \brk<c_{1}, B_{1}>]$ and 
	$\mathbbm I_{2}=[\brk<\one, \Set{\one}>, \brk<c_{2}, B_{2}>]$ are 
	order-isomorphic iff $B_{1}$ and $B_{2}$ are isomorphic. 
\end{lem}
\begin{proof}
	Suppose that $\mathbbm I_{1}\simeq \mathbbm I_{2}$. Let $A_{i}$ be 
	the MR-subalgebra located by $\brk<c_{i}, B_{i}>$. Let
	$a_{i}$ be an $A_{i}$-vertex and let $1>s_{1}>\dots>s_{k_{i}}=a_{i}$ 
	be a maximal chain in $A_{i}\cap[a_{i}, \one]$. Then we have  a 
	maximal chain of subalgebras of $A_{i}$ induced by the intervals 
	$[s_{j}, \one]\cap A_{i}$ -- this is a maximal chain as the 
	$A_{i}$-rank goes up by one as $j$ increases by one. 
	
	From this we see that the rank of $A_{1}$ is equal to that of $A_{2}$. 
	As the rank of $B_{i}$ equals that of $A_{i}$ we have $B_{1}\simeq 
	B_{2}$. 
	
	Conversely, if $B_{1}\simeq B_{2}$ then we have 
	$A_{1}\simeq \rsf I(B_{1})\simeq \rsf I(B_{2})\simeq A_{2}$ from 
	which the result is clear. 
\end{proof}

Now we need to count the number of MR-algebras with locator of a 
particular dimension. 

\begin{lem}
    There are $2^{n-m}S(n, m)$ MR-subalgebras of $\mathcal L$ with 
    dimension $m$ and locator a subalgebra of $[0,1]$. 
\end{lem}
\begin{rem}
    $S(n, m)$ is a Stirling number of the second kind,  counting the 
    number of partitions of a set of size $n$ into $m$ pieces.
\end{rem}
\begin{proof}
    There are $S(n, m)$ Boolean subalgebras of $[0,1]$ of dimension 
    $m$ -- since each subalgebra corresponds to a partition of the 
    atoms of $[0,1]$ into $m$ pieces. 
    
    Given such a subalgebra, we see that $\brk<c,B>\sim\brk<c', B>$ 
    iff $c+c'\in B$ iff $c$ and $c'$ are in the same coset in $[0,1]$ 
    relative to $B$.  Thus the number of cosets of $B$ equals the 
    number of MR-subalgebras located by $B$, ie $2^{n-m}$. 
    
    Hence there are $2^{n-m}S(n, m)$ such MR-subalgebras. 
\end{proof}

Now we are able to compute the M\"obius function. First a small lemma. 
\begin{lem}\label{lem:cloMR}
	Let $\mathcal M$ be the partial order of MR-subalgebras of $\mathcal 
	L$. Then 
	$\comp{\mathcal M}$ is isomorphic to the partial order of 
	implication sublattices of $[0,1]$. 
\end{lem}
\begin{proof}
	We know that $\comp{\brk<c, B>}=\brk<c, B>$ iff $\brk<c, 
	B>\sim\brk<\min B, B>$.  Thus the mapping that takes a closed 
	element to the second component of a locator pair is an order 
	isomorphism. 
\end{proof}

\begin{cor}\label{cor:cloMR}
	$$
	\mu_{\comp{\mathcal M}}(\Set{\one}, \mathcal L)=n!. 
	$$
\end{cor}
\begin{proof}
	Noting that $\Set{\one}$ is closed we can apply the lemma and 
	equation
	\eqref{eq:muOne}. 
\end{proof}

Using this result and \thmref{thm:closure} we see that 
\begin{align*}
    \sum_{\comp{\brk<c,B>}=\mathcal L}\mu(\Set{\one}, \brk<c,B>)&=n!\\
    \intertext{ and so}
    \mu(\Set{\one}, \mathcal L_{n}) &= n!- 
    \sum_{\substack{\comp{\brk<c,B>}=\mathcal L\\
    B\not=[0,1]}}\mu(\Set{\one}, \brk<c,B>)\\
    &= n!- 
    \sum_{m=\lceil\frac{n}2\rceil}^{n-1}\sum_{\substack{\brk<c,B>\\ 
    \dim B=m}}\mu(\Set{\one}, \brk<c,B>)\\
    &= n!- 
    \sum_{m=\lceil\frac{n}2\rceil}^{n-1}2^{n-m}S(n, m)\mu(\Set{\one}, \mathcal L_{m}). 
\end{align*}

Now let $a_{n}=\mu(\Set{\one}, \mathcal L_{n})$ so we can rewrite this 
as 
$$
\frac{a_{n}}{2^{n}n!}=\frac1{2^{n}}-\sum_{m=\lceil\frac{n}2\rceil}^{n-1}(m! S(n, m))\frac{a_{m}}{2^{m}m!}
$$
or as 
$$
\sum_{m=\lceil\frac{n}2\rceil}^{n}S(n, m)\frac{a_{m}}{2^{m}}= n!. 
$$

\end{document}